\documentclass[11pt]{article}

\usepackage{amsmath}
\usepackage{amssymb}
\usepackage{amsfonts}
\usepackage{amsthm}

\usepackage[utf8x]{inputenc} 

\usepackage{graphicx}
\usepackage[autostyle]{csquotes}
\usepackage[linktocpage]{hyperref}
\usepackage{fullpage}
\usepackage{slashed}
\usepackage{mathtools}
\usepackage{indentfirst}
\usepackage{bm}
\usepackage{mathrsfs} 
\usepackage{authblk}
\usepackage{tikz}

\usepackage{subfig}

\allowdisplaybreaks

\newcommand{\intg}{\mathbb{Z}}

\newcommand{\real}{\mathbb{R}}
\newcommand{\complex}{\mathbb{C}}

\newcommand{\rarw}{\rightarrow}

\newcommand{\imp}{\Rightarrow}

\newcommand{\lb}{\left(}
\newcommand{\rb}{\right)}
\newcommand{\lsb}{\left[}
\newcommand{\rsb}{\right]}

\newcommand{\RN}[1]{%
  \textup{\uppercase\expandafter{\romannumeral#1}}%
}

\title{\textbf{Witt invariants from $q$-series} $\hat{Z}$}
\author{John Chae}
\affil{Center for Quantum Mathematics and Physics (QMAP), UC Davis, Davis, CA, 95616, USA \\ {yjchae@ucdavis.edu}}
\date{}  

\begin{document}

\maketitle

\begin{abstract}
We present a relation between the Witt invariants of 3-manifolds and the $\hat{Z}$-invariants. It provides an alternative approach to compute the Witt invariants of 3-manifolds, which were originally defined geometrically in four dimensions. We analyze various homology spheres including a hyperbolic manifold using this method.
\end{abstract}

\tableofcontents

\section{Introduction}

In the last few years, there has been an extensive development on the $q$-series named $\hat{Z}$ introduced in \cite{GPV, GPPV}, which are conjectured to categorify the 
Witten-Reshitikhin-Turaev (WRT) invariant~\cite{W, RT} for a closed oriented 3-manifold $Y$. It was shown that $\hat{Z}$ invariant have multifaceted characteristics across low dimensional topology, number theory and mathematical physics. Consequently, $\hat{Z}$ established links between these three fields. From mathematical physics perspective, an existence of $\hat{Z}$ is predicted by a 3-dimensional $\mathcal{N}=2$ supersymmetric quantum field theory, which arose from a compactification of M-theory~\cite{GPV, GPPV}. More precisely, the structure of the field theory's  BPS sector of the Hilbert space
$$
\mathcal{H}_{BPS}(Y)\quad = \bigoplus_{b,i,j} \mathcal{H}^{i,j}_{BPS}(Y;b),
$$ 
predicts a $q$-series with integer coefficients as the graded Euler characteristic of $\mathcal{H}_{BPS}(Y)$
$$
\chi[\mathcal{H}^{i,j}_{BPS}(Y;b)]=\hat{Z}_{b}[Y;q]= \sum_{i,j} (-1)^i\, q^j\, \text{dim}\, \mathcal{H}^{i,j}_{BPS}(Y;b),\qquad b \in Spin^c (Y).
$$
The integer coefficients reflect the number of BPS states. Furthermore, $\hat{Z}_{b}$ was developed with having in mind a long-term goal of a categorification of the $SU(2)$ WRT invariant~\footnote{The normalization used here is $\tau [S^3;k]=1$.}
\begin{equation}
\tau [Y;k]= \frac{1}{i\sqrt{2k}} \sum_{a,b \in Spin^c (Y)/ \intg_2 } c^{WRT}_{ab} \, \hat{Z}_{b}(q) \Big\vert_{q \rarw e^{\frac{i2\pi}{k}}}.
\end{equation}
\newline
The $q$-series enjoys the $Spin^c$ conjugation symmetry
$$
\hat{Z}_{-b}(q) = \hat{Z}_{b}(q).
$$
\indent The number theoretic aspect of $\hat{Z}$ is manifested by its modular property. In \cite{CCFGH, HJC}, it was demonstrated that $\hat{Z}$ for Siefert fibered rational homology spheres can be expressed in terms of false or mock theta functions of certain weights, which are well-known examples of quantum modular forms (see \cite{Z} for a review). This modular characteristic of $\hat{Z}$ has uncovered the origin of the modular property of the WRT invariant, which was discovered for the Poincare homology sphere in \cite{LZ} and then generalized to an arbitrary integral homology sphere in \cite{H2}. From topology viewpoint, $\hat{Z}$ inspired an introduction of an invariant of a knot complement $M_K$ in \cite{GM}, which is a two variable series:
$$
F_K (x,q) := \hat{Z}[M_K] \in 2^{-c} q^{\Delta} \intg [ x^{\pm 1/2}] [[q^{\pm 1}]]
$$
This series has broadened the range of 3-manifolds for which $\hat{Z}$ can be computed through Dehn surgery, including hyperbolic 3-manifolds. 
Later, it was shown in \cite{GHNPPS} that this series invariant $F_K$ in turn are connected to the Akutsu-Deguchi-Ohtsuki (ADO) polynomials~\cite{ADO}. This conjecture was reinforced in \cite{C}. Recently, another facet of $\hat{Z}$ was revealed in \cite{GPP}. Specifically, the authors of \cite{GPP} investigated the spin refined version of the WRT invariant at the fourth root of unity and elucidated that the corresponding $\hat{Z}$'s are related to the Rokhlin invariant $\mu(Y,s)$ and the d-invariant (or the correction term) of a certain version of the Heegaard Floer homology for several classes of 3-manifolds (see Figure 1),
$$
e^{-i 2\pi \frac{3\mu(Y,s)}{16}} = \sum_{b} c^{Rokhlin}_{sb} \hat{Z}_{b}[Y;q] \bigg|_{q \rarw i}\qquad s \in Spin(Y),\quad b \in  Spin^c (Y).
$$
This work has exemplified an existence of a connection between geometric topology and the quantum invariant $\hat{Z}$.
\newline

In this paper, we find a new relation of the same type of connection. Namely, a link between the Witt invariant $w(Y)$, Witt defect $\text{def}_{3}(\Theta)$ and $\hat{Z}$ of $Y$ from a certain refinement of the WRT invariant at the sixth root of unity (see Section 2 for a review). The two former invariants are geometrically defined on the level of 4-manifolds, thus they also posses cobordism characteristic.
$$
 i^{-w(Y) + 2\Theta^3 +  \text{def}_{3}(\Theta)} \sqrt{3}^{\epsilon(\Theta) + d(Y_{\Theta}) - d(Y)} =  \sum_{b\in Spin^c (Y) / \intg_2} c^{Witt}_{\Theta b} \hat{Z}_{b} (q)\bigg|_{q \rarw e^{\frac{i2\pi}{6}}}\qquad \Theta \in H^{1}(Y; \intg / 2\intg).
$$
For rational homology spheres $Y$ ($H_1 (Y; \intg) = \intg / p\intg $), there are two different cases. The first case is when $p=$ odd, 
\begin{equation}
c^{Witt}_{t}= \frac{e^{-i\pi/4}}{4\sqrt{3}} \sum_{r=0}^{5} e^{-\frac{i\pi}{12p} (2pr-2t+p)^2},
\end{equation}
where $t=0,\cdots, p-1$. When $p=$even,
\begin{equation}
c^{Witt}_{w t}= \frac{e^{-i\pi/4}}{2\sqrt{3}}e^{-\frac{i\pi}{3p} (t+\frac{p}{2}(w+1))^2} \lb 1 + e^{\frac{i\pi}{3} (pw+2t)} + e^{\frac{i2\pi}{3} (pw+2t-p)} \rb ,
\end{equation}
where $w=0,1$ and $t=0,\cdots, p-1$~\footnote{We used the fact that $Spin^c (Y)$ is affinely isomorphic to $H_1 (Y; \intg)$}. This new relation not only enriches the conceptual aspects of the invariants, it provides a new method of computing the Witt invariant and Witt defect directly in three dimension as well.
\begin{figure}[h]
\centering
\includegraphics[scale=0.5]{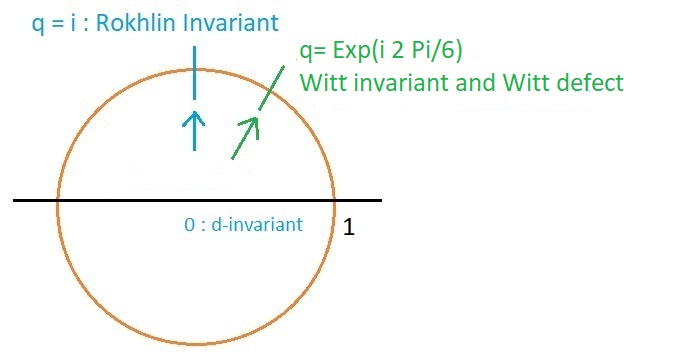}
\caption{Topological invariants at the fourth and the sixth roots of unity from the limits of $\hat{Z}$.}
\end{figure}
\newline

The rest of the paper is organized as follows.
In Section 2, we give a review of the Witt invariant, and Witt defect. Then, in Section 3, we express the refined WRT invariant at six root of unity in terms of $\hat{Z}$ for plumbed 3-manifolds.  In Section 4, we apply the new relation to several classes of 3-manifolds to obtain their Witt invariant and Witt defect.


\section{The Witt invariant and defect}

Let us review the Witt invariant and defect of 3-manifolds defined in \cite{KMZ}. Their formulation takes place in four dimension. Let $Y$ be a closed oriented 3-manifold. By the vanishing of its oriented cobordism group $\Omega(Y)=0$~\cite{GS}, $Y$ bounds a compact oriented 4-manifold $X$ whose intersection form is denoted by $\tilde{Q}_X$. Its signature is denoted by $\sigma(X)$. We next diagonalize $\tilde{Q}_X$ in $\intg_3$-coefficient ring, obtaining $0, \pm 1$ as its diagonal entries. We denote it by $Q_X$. Then we let $w(X)$ to be its trace Tr $Q_X$. The mod 3 Witt invariant of $Y$ is defined as
\begin{equation}
w(Y) : = \sigma(X) -  w(X)\quad \text{mod}\, 4.
\end{equation}
$w(Y)$ is independent of $X$. Since we deal with a compact 4-manifold with a boundary, we would like to detect an effect of the boundary. This leads to the notion of the Witt defect. Specifically, we consider a cyclic n-fold cover manifold $\tilde{Y} \rarw Y$. By the result of \cite{CG}, this covering manifold extends to a cyclic branched cover $\tilde{X} \rarw X$ branched along a closed surface $F$ in $X$. We let $Q_{\tilde{X}}$ be an intersection form of $\tilde{X}$ in $\intg_3$ coefficient. The mod 3 Witt defect of $\tilde{Y} \rarw Y$ is defined as
$$
\text{def}_{3}(\tilde{Y} \rarw Y) := n w(X) - w(\tilde{X}) - \frac{n^2 -1 }{3n} F \cdot F \quad \text{mod}\, 4,
$$
where n divides $F \cdot F$. The specific Witt defect that is relevant in our context is a double cover 3-manifold equipped with a cohomological class $\Theta \in H^{1}(Y ; \intg/2\intg)$:
\begin{equation}
\text{def}_{3}(Y_{\Theta} \rarw Y) = 2 w(X) - w(\tilde{X}) - \frac{1}{2} F \cdot F \quad \text{mod}\, 4.
\end{equation}
We abbreviate the above defect as $\text{def}_{3}(\Theta)$. Due to the presence of the boundary, the difference between the first two terms in (5) is not necessarily zero. Note that $w(Y)$ and $\text{def}_{3}(Y_{\Theta} \rarw Y)$ taking value in $\intg/4\intg$ follows from the fact that the Witt ring $W(R)$ of $R=\intg/3\intg$ is $\intg/4\intg$~\cite{MH}.

\section{A refined WRT Invariant at the sixth root of unity}

In this section we derive a relation between the $H^{1}$-refined $SU(2)$ WRT invariant~\cite{KM} and the $\hat{Z}$-invariant at the sixth root of unity $k= 2 $ mod 4 for a negative definite plumbed 3-manifold $Y$. Specifically, this invariant deals with a closed oriented 3-manifold equipped with a 1-dimensional cohomological class $w \in H^{1}(Y; \intg / 2\intg)$. Let us review plumbed 3-manifolds and the corresponding refined WRT invariant at an even root of unity ($k=$even).
\newline

Plumbed 3-manifolds are characterized by a weighted graph $\Gamma$. Each vertex of the graph represent a $S^1$-bundle over a compact g-surface and is labeled by $[g,n \in \intg]$, where $n$ is the Euler number of the bundle. An edge corresponds to a gluing between two $S^1$-bundles in the fiber-base exchanging way. We choose all the base surfaces to be $S^2$'s ($g=0$ suppressed from now on). From $\Gamma$, a link $L(\Gamma)$ can be obtained by replacing each vertex of $\Gamma$ by an unknot whose framing $n_i$ is set by the Euler number and each edge of $\Gamma$ by a Hopf link between two vertices (Figure 2). Performing a Dehn surgery on $L(\Gamma)$ yields $Y$.
\begin{figure}[h]
\centering
\includegraphics[scale=0.5]{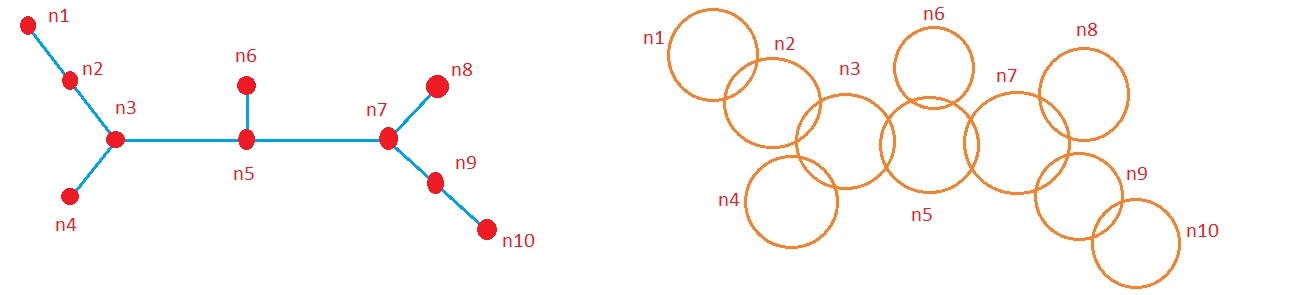}
\caption{A plumbing graph $\Gamma$ and a link $L(\Gamma)$ associated to $\Gamma$.}
\end{figure}
\noindent Furthermore, the adjacency matrix $B$ of $\Gamma$ coincides with the linking matrix of $L(\Gamma)$. Its entries are given by
$$
B_{ij} = \begin{cases}
n_i,\quad i=j\\
1,\quad i\, \&\, j\, \text{are connected}\qquad i,j \in V\\
0,\quad \text{otherwise}
\end{cases}
$$
where $V$ is a set of vertices of $\Gamma$. In this paper we focus on $Y(\Gamma)$ whose $H_1 (Y(\Gamma) ; \intg) = \intg^L / B \intg^L $, where $L$ is the number of components of $L(\Gamma)$\footnote{Since it is clear from the context whether $L$ stands for the link $L(\Gamma)$ or the number of link components, we use $L$ interchangeably.} 
\newline

The $H^{1}$-refined WRT invariant of $Y=Y(\Gamma)$ for $k=$ even is given by~\cite{GPP}
$$
\tau_k [Y(\Gamma),k,s_c] = \frac{F^{(c+ \epsilon)}[L(\Gamma)]}{F[U_{+1}]^{b_{+}} F[U_{-1}]^{b_{-}}},\qquad  s_c \in H^{1}(Y(\Gamma); \intg / 2\intg),
$$
where $U_{\pm 1}$ are $\pm 1$-framed unknots, $b_{\pm}$ are the number of positive and negative eigenvalues of $B(\Gamma)$ and 
\begin{align*}
F^{(c+ \epsilon)}[L(\Gamma)]  & = \sum^{\prime}_{\substack{ 0 \leq n_{r} \leq 2k-1 \\ n_r = c_r } } J_{n}(L(\Gamma); q)  \prod_{s=1}^{L} \frac{q^{n_{s}/2}- q^{-n_{s}/2}}{q^{1/2}- q^{-1/2}} \bigg|_{q=e^{i 2\pi/k}}\qquad c \in \intg^L ,\quad \epsilon = (1,1,\cdots,1) \in \intg^L\\
F[L(\Gamma)] & =  \sum^{\prime}_{ 0 \leq n_{r} \leq 2k-1 }  J_{n}(L(\Gamma); q)  \prod_{s=1}^{L} \frac{q^{n_{s}/2}- q^{-n_{s}/2}}{q^{1/2}- q^{-1/2}}  \bigg|_{q=e^{i 2\pi/k}},\qquad n \in \intg_{+}^L
\end{align*}
The prime in the summations means $n_r = 0, k$ are omitted (they result in diverging  $\tau_k$). $J_{n}(L(\Gamma); q)$ is the $SU(2)$ colored Jones polynomial of $L(\Gamma)$:
$$
J_{n}(L(\Gamma); q) = \frac{1}{q^{1/2}-q^{-1/2}} \prod_{r \in V} q^{\frac{a_r (n_{r}^{2} -1 )}{4}} \lb \frac{1}{q^{n_r /2}-q^{- n_r /2}} \rb^{\text{deg}(r) - 1} \times\\
\prod_{(i,j) \in \text{E}} \lb q^{n_i n_j /2} - q^{- n_i n_j /2} \rb,
$$
where $E$ is a set of edges of $\Gamma$.
\newline

In order to express $\tau_k$ in terms of $\hat{Z}$, the strategy from Appendix A of \cite{GPPV} is applied. The core parts of the strategy are the analytic continuation of $q$-domain ($|q| = 1$) to the complex unit disc ($ |q| \leq 1 $) and the application of the Gauss sum reciprocity formula (\cite{GPP} (2.25))\footnote{The refined WRT invariant requires a different version of the reciprocity formula than the one used in \cite{GPPV}}. After this procedure, we arrive at (2.26) of \cite{GPP}:
$$
\tau [Y, w_c ;k]= \frac{1}{Sin(\pi/k) 2^{L+2}|\text{Det}\, B|^{1/2}} \sum_{\tilde{a},\tilde{b} \in \intg^{L}/2B\intg^{L}} e^{-i\pi \frac{k}{2} \tilde{a}^T B^{-1}\tilde{a}} e^{-i\pi \tilde{a}^T B^{-1}\tilde{b}} \hat{Z}_{\tilde{b} + B(c+\epsilon)}\bigg|_{q \rarw e^{\frac{i2\pi}{k}}},
$$
where
$$
\hat{Z}_b [Y(\Gamma) ; q] = (-1)^{b_{+}} q^{\frac{3\sigma - Tr B}{4}} \sum_{w \in 2B\intg^L + b} F_{w}\, q^{-\frac{w^T B^{-1} w}{4}}\, \in 2^{-m}q^{\Delta(b)}\intg[[q]],\qquad b \in Spin^c (Y)
$$
$$
\sum_{w \in 2B\intg^L + \delta} F_{w} \prod_{r \in V} x_{r}^{w_r} := \prod_{r \in V} \lsb  \frac{1}{\lb x_r - x_{r}^{-1} \rb^{deg(r) -2}} +  \frac{1}{\lb x_r - x_{r}^{-1} \rb^{deg(r) -2}} \rsb,\qquad \delta(i) = deg(i)\quad \text{mod}\, 2.
$$
The first and the second term in the above square bracket are expanded around $x=0$ and $x=\infty$, respectively. For our purpose, we set $k=2$ mod 4 from here and split $\tilde{a}=a +BA,\, a \in \intg^L / B \intg^L,\, A \in \intg_{2}^L $. Then,
$$
\tau [Y,w_c ;k]= \frac{-i}{Sin(\pi/k) 2^{L+2}|\text{Det}\, B|^{1/2}} \sum_{\substack{ a\in \intg^{L}/B\intg^{L}\\  \tilde{b}\in \intg^{L}/2B\intg^{L}\\ A \in \intg_{2}^{L}}} e^{-i\pi \frac{k}{2} a^T B^{-1}a} e^{-i\pi a^T B^{-1}\tilde{b}} e^{-i\pi A^T B A -i\pi A^T\tilde{b}}\, \hat{Z}_{\tilde{b} + B(c+\epsilon)}\bigg|_{q \rarw e^{\frac{i2\pi}{k}}}.
$$
We next apply $\sum_{j} B_{ij}c_{j}= 0$ mod 2 and $\delta_{i}=\sum_{j\neq i}B_{ij}$, where $\delta_{i} = deg(i)$ mod 2. This yields
$$
\tau [Y,w_c ;k]= \frac{-i}{4 Sin(\pi/k) |\text{Det}\, B|^{1/2}} \sum_{\substack{ a\in \intg^{L}/B\intg^{L}\\  \tilde{b}\in \intg^{L}/2B\intg^{L}}} e^{-i\pi \frac{k}{2} a^T B^{-1}a} e^{-i\pi a^T B^{-1}\tilde{b}}\, \delta ( \tilde{b} =\beta \quad \text{mod}\quad 2) \, \hat{Z}_{\tilde{b} +B(c+\epsilon)}\bigg|_{q \rarw e^{\frac{i2\pi}{k}}},
$$
$$
\beta = \begin{cases}
0,\quad | \text{Det}\, B | = \text{even} \\ 
1,\quad | \text{Det}\, B | = \text{odd} \\ 
\end{cases}
$$
In the case of $\beta = 0 $, we let $\tilde{b} = 2b$,
$$
\tau [Y,w_c ;k]= \frac{-i}{4 Sin(\pi/k) |\text{Det}\, B|^{1/2}} \sum_{\substack{ a\in \intg^{L}/B\intg^{L}\\ b\in \intg^{L}/ B\intg^{L}}} e^{-i\pi \frac{k}{2} a^T B^{-1}a} e^{-i 2\pi a^T B^{-1}b}\, \hat{Z}_{2b +B(c+\epsilon)}\bigg|_{q \rarw e^{\frac{i2\pi}{k}}}.
$$
When $\beta = 1 $, we split $\tilde{b} = 2b + \epsilon $,
$$
\tau [Y,w_c ;k]= \frac{-i}{4 Sin(\pi/k) |\text{Det}\, B|^{1/2}} \sum_{\substack{ a\in \intg^{L}/B\intg^{L}\\  b \in \intg^{L}/ B\intg^{L}}} e^{-i\pi \frac{k}{2} a^T B^{-1}a} e^{-i2\pi a^T B^{-1}b} e^{-i\pi a^T B^{-1} \epsilon}\, \hat{Z}_{2b + \epsilon +B (c+\epsilon)}\bigg|_{q \rarw e^{\frac{i2\pi}{k}}}.
$$
We propose that for any rational homology sphere $Y$ and $k=2$ mod 4,
\begin{equation}
\tau [Y,w; k]= \frac{-i}{4 Sin(\pi/k) |H_{1}(Y; \intg)|^{1/2}} \sum_{a,b \in H_{1}(Y; \intg)} e^{-i2\pi \frac{k}{4} lk(a,a)} e^{-i 2\pi lk(a, b + \frac{\lambda}{2})}\,  \hat{Z}_{\phi (w,b)}\bigg|_{q \rarw e^{\frac{i2\pi}{k}}},
\end{equation}
$$
\lambda = \begin{cases}
0, \quad | H_{1}(Y; \intg) | = \text{even} \\ 
\epsilon ,\quad | H_{1}(Y; \intg) | = \text{odd} \\ 
\end{cases}
$$
$$
\phi  : H^{1}(Y;\intg / 2\intg) \times H_{1}(Y; \intg) \rarw Spin^{c}(Y)
$$
This $\phi$ map converts a 1-cocycle in $H^{1}(Y;\intg / 2\intg)$ to a $Spin^{c}$ structure of $Y$.

\section{Witt Invariants and the $q$-series}

In this section, we relate the Witt invariant and the Witt defect to $\hat{Z}$ and find sum rules. In \cite{KMZ}, the $SU(2)$ WRT invariant at the sixth root of unity for a closed oriented 3-manifold was investigated. It was shown that the WRT invariant is a sum of the invariants of the manifold equipped with a 1-dimensional mod 2 cohomology class $\Theta$:
$$
\tau_{6}[Y] = \sum_{\Theta \in H^{1}(Y; \intg/ 2\intg)} \tau_{6}[Y, \Theta].
$$
Furthermore, $\tau_{6}[Y, \Theta]$ can be expressed in terms of $w(Y)$ and $\text{def}_{3}(\Theta)$,
\begin{equation}
\tau_{6}[Y, \Theta] = i^{-w(Y) + 2\Theta^3 + \text{def}_{3}(\Theta)} \sqrt{3}^{\epsilon(\Theta) + d(Y_{\Theta}) - d(Y)}
\end{equation}
\newline
$$
d(Y) = rk H^{1}(Y; \intg/ 3\intg), \qquad  d(Y_{\Theta}) = rk H^{1}(Y_{\Theta}; \intg/ 3\intg), \qquad 2\Theta^3 \in \intg / 4\intg
$$
$$
w(Y): \text{mod 3 Witt invariant of}\, Y \quad (cf. (4)) 
$$
$$
\text{def}_{3}(\Theta) : \text{mod 3 Witt defect of the double cover manifold}\, Y_{\Theta} \rarw Y \quad (cf. (5))
$$
$$
\epsilon(\Theta)=   \begin{cases}
0,\quad \Theta =0 \\
1,\quad \Theta \neq 0\\
\end{cases}
$$
\newline
\noindent The $\hat{Z}$ decomposition of $\tau_{6}$ is
\begin{equation}
\tau_{6}[Y] = \sum_{\Theta \in H^{1}(Y; \intg/ 2\intg)} \sum_{b\in Spin^c(Y)/ \intg_2} c^{Witt}_{\Theta b} \hat{Z}_{b}(q)\bigg|_{q \rarw e^{\frac{i2\pi}{6}}}.
\end{equation}
and let
\begin{equation}
\tau_{6}[Y, \Theta] : =  \sum_{b\in Spin^c(Y)/ \intg_2} c^{Witt}_{\Theta b} \hat{Z}_{b} (q)\bigg|_{q \rarw e^{\frac{i2\pi}{6}}} .
\end{equation}
For generality, we utilize the unfolded versions of (8) and (1) to deduce a consistency condition as sum rules. 
\begin{equation}
\sum_{\Theta \in H^{1}(Y ; \intg /2\intg)} c^{Witt}_{\Theta t} = \sum_{a \in Spin^c (Y)} c^{WRT}_{a t} (k=6) \equiv c^{WRT}_{t} (k=6).
\end{equation}
The Witt coefficients in (2) and (3) must satisfy (10). For a rational homology sphere $H_1 (Y; \intg) = \intg / p\intg $,\, $c^{WRT}_{a t}(k)$ is given by~\cite{GPP}
$$
c^{WRT}_{a t}(k) = e^{-i 2\pi k \, \text{lk}(a,a)} \frac{e^{\frac{i 4\pi at}{p}}}{\sqrt{p}},\qquad a,t=0, \cdots , p-1
$$
where $\text{lk}(a,a)$ is the linking form on Tor $H_1 (Y; \intg)$.

%

\section{Examples}

In this section, we compute the Witt invariant, Witt defect and other invariants for homology spheres using the $q$-series $\hat{Z}$.

\subsection{Brieskorn spheres}

We begin our analysis with the Brieskorn spheres. Since they are $\intg HS^3$, they have a unique $Spin^c$ structure. First, we consider the following family 
$$
Y_r = S^{3}_{-\frac{1}{r}}(3^{r}_{1}) = \Sigma(2,3,6r+1), \qquad r \in \intg_{+}
$$
where $3^{r}_{1}$ is the right handed trefoil. These manifolds have $H_1 (Y;\intg)=0$ so they carry a unique $Spin^c$ structure. Their $\hat{Z}$ in terms of the (quantum) modular forms can be computed using the general formula in \cite{GM} (Proposition 4.8).
$$
\hat{Z} [Y_r ;q] = q^{\frac{1}{2} - \frac{(6r-5)^2}{24(6r+1)}} \lb \Psi^{(6r-5)}_{6(6r+1)} (q) - \Psi^{(6r+7)}_{6(6r+1)}(q) - \Psi^{(30r-1)}_{6(6r+1)}(q) + \Psi^{(30r+11)}_{6(6r+1)}(q) \rb.
$$
Using (9) together with $c^{Witt}_{00}(6) = -i/2$, we obtain 
$$
\tau_{6}[Y_r] = \tau_{6}[Y_r, 0] = 1.
$$
The desired invariants can be read off from (7),
\begin{equation}
w(Y)=0\quad \text{mod},\, 4 \qquad d(Y)=d(Y_0) = 0.
\end{equation}
The vanishing of the $d$'s is clear since $Y=\intg HS$ and $d(Y_0) = 2d(Y)$. Since (4) is cobordism in nature, meaning of $w(Y)$-value is tied to behaviors of bounding $X$ and its cover $\tilde{X}$. Specifically, from (4), we get $\sigma(X) \equiv w(X)$ mod 4. Using the fact that any $M=\intg HS$ carries a unique spin structure $s_0$ and bounds a smooth spin 4-manifold $W$\footnote{$\intg_2$ homology spheres can bound a smooth $W$~\cite{S}. Moreover, the vanishing of the spin cobordism group of 3-manifolds $\Omega^{spin}_3 = 0$~\cite{GS} implies that a $Spin(M)$ structure extends to a $Spin(W)$ structure.}, an application of the Rokhlin's theorem~\cite{R} $\mu (M,s) \equiv \sigma(W)\ \text{mod}\, 16,\,  s\in Spin(M)$, we arrive at $\mu (Y,s_0) = w(X)$ mod 4. So $w(X)$ determines whether $Y$ admits a smooth embedding in $S^4$. We next apply the generalization of the (Hirzebruch) signature theorem in \cite{APS1}, $\sigma(X) = p_1 (TX)/3 - \eta(0)$, where $p_1 \in H^4(TX)$ is the first Pontryagin class of the tangent bundle $TX$ and $\eta(s)$ is a spectral invariant\footnote{$\eta(s)$ is defined by eigenvalues of a first order differential operator on $W$ (\cite{APS1} (1.7)). Note that without $\eta(0)$ term, it is the classical Hirzebruch signature theorem~\cite{MS}.} measuring an effect of the boundary $Y$. From the theorem, we get $p_1(TX) = 3w(X) + 3\eta(0)$ mod 12. By another theorem of Rokhlin~\cite{MK}, $p_1 (TX) \in 48 \intg $, which constrains $w(X)$ and/or $\eta(0)$. In case $p_1(TX)=0$, hence $p_1(TX)/2=0$ implies that $X$ is a string manifold~\cite{DHH}. So  $w(X)$ together with $\eta(0)$ controls whether or not $X$ can carry a string structure. 
\newline

\noindent We next focus on
$$
Y_r = S^{3}_{-\frac{1}{r}}(3^{l}_{1}) = \Sigma(2,3,6r-1), \qquad r \in \intg_{+}
$$
where $3^{l}_{1}$ is the left handed trefoil. As in the previous example, $\hat{Z}$ for this family of manifold are calculated as~\cite{GM}.
$$
\hat{Z} [Y_r;q] = q^{-\frac{1}{2} - \frac{(6r+5)^2}{24(6r-1)}} \lb \Psi^{(6r+5)}_{6(6r-1)} (q) - \Psi^{(6r+7)}_{6(6r-1)}(q) - \Psi^{(30r-11)}_{6(6r-1)}(q) + \Psi^{(30r+1)}_{6(6r-1)}(q) \rb.
$$
Using (7) and (9), we arrive at the same result as that of $\Sigma(2,3,6r+1)$. Therefore, the Witt invariant and defect are insensitive to the chirality of these surgery knots. We now choose other torus knots. Specifically, we pick 
\begin{equation}
S^{3}_{-\frac{1}{r}}(T(2,\pm 5)), \qquad S^{3}_{-\frac{1}{r}}(T(2,\pm 7)), \qquad S^{3}_{-\frac{1}{r}}(T(3,\pm 5)), \qquad S^{3}_{-\frac{1}{r}}(T(3,\pm 7)).
\end{equation}
Their $\hat{Z}$ can be found in the same way as the above examples (see Appendix B for their explicit expressions). In all cases, their $w(Y)$, $d(Y)$ and $d(Y_0)$ are same as (11). This result can also be deduced from the fact that the manifolds in (12) are all $\intg HS$ and $\tau_{6}[\intg HS]=1$, which was shown in \cite{KMZ}~\footnote{Since $H^{1}(\intg HS ;\intg_3)=0$ so $d(Y)=0$, which sets the above value of $w(Y)$.}.
\newline

We finish this subsection with the Poincare homology sphere $P=\Sigma(2,3,5)$. It can be obtained by plumbing on $-E_8$ graph.
Its $\hat{Z}$ is given~\cite{GM}
$$
\hat{Z}[\Sigma(2,3,5);q]= q^{-3/2} \lb 2- \sum_{n=1}^{\infty} \chi_{+}(n) q^{(n^2-1)/120} \rb,
$$
where 
$$
\chi_{+}(n) := \begin{cases}
1,\quad n \equiv 1,11,19,29 \quad \text{mod}\, 60  \\
-1,\quad n \equiv 31,41,49,59 \quad \text{mod}\, 60\\
0,\quad \text{otherwise}
\end{cases}
$$
The Witt invariant $w(P)$ and $d(P)$ coincide with (11) for the same reason as (12) (i.e. $P=\intg HS^3$).

\subsection{Lens spaces}

A well known and the simplest rational homology sphere is a Lens space $Y= -L(p,1)\, (p > 1)$. They carry $|H_{1} (Y;\intg)|$ number of $Spin^c$ structures. Since these manifolds are described by one vertex plumbing graph, it is straightforward to obtain their nonzero $\hat{Z}_b$ using (15) in Appendix C.
$$
\hat{Z}_0 [Y;q] = -2 q^{\frac{p-3}{4}},\qquad \hat{Z}_1 = 2 q^{\frac{p^2 -3p+4}{4p}}.
$$
We next compute $H^{1}(Y; \intg/ 2\intg)$ via the universal coefficient theorem for cohomology
\begin{equation*}
\begin{aligned}
H^{1}(Y;\intg/ 2\intg) \cong &\, \text{Hom} ( H_{1}(Y; \intg) ,\intg/ 2\intg) \oplus \text{Ext} ( H_{0}(Y; \intg), \intg/ 2\intg)\\
\cong & \, \intg/ p\intg \otimes \intg/ 2\intg  \\
\cong & \, \intg / \text{gcd}(2,p) \intg
\end{aligned}
\end{equation*}
where $\text{Ext} ( \intg ,  \intg/ 2\intg) =0$ is used. 
When $p=$ odd, $H^{1}(Y;\intg/ 2\intg) = 0 $ so (9) becomes
$$
\tau_{6}[Y] = \tau_{6}[Y, \Theta = 0 ] =  \sum_{b\in Spin^c(Y) /\intg_2} c^{Witt}_{b} \hat{Z}_{b} (q = e^{\frac{i2\pi}{6}}).
$$
After straightforward calculation we obtain $\tau_{6}[Y, \Theta = 0 ]$ for -L(p,1), p=3,5,7, which are given in the following table
\begin{center}
\begin{tabular}{ |c|c| } 
\hline
$-L(p ,1 )$       & $ \tau_{6}[\Theta=0]$  \\ 
 \hline
$-L(3,1)$ & $i\sqrt{3}$  \\ 
 \hline
$-L(5,1)$ & $-1$ \\ 
 \hline
$-L(7,1)$ & $ 1 $ \\ 
\hline
\end{tabular}
\end{center}
This table agrees with the results of \cite{KMZ}\footnote{Our orientation convention for the manifold is opposite of that of \cite{KMZ}}. This WRT invariant $\tau_{6}$ has a period from $p=3$ to $p=7$; for instance,
$$
\tau_{6}(p=9) = \tau_{6}(p=3),\quad \tau_{6}(p=11) = \tau_{6}(p=5)\quad \text{etc}.
$$
We next use (7) to find $w(Y)$ and $d(Y)$:
\begin{center}
\begin{tabular}{ |c|c|c| } 
\hline
$-L(p ,1 )$  & $w(Y) \in \intg /4\intg $ & $d(Y) \in \intg $ \\ 
 \hline
$-L(3,1)$    & 3   & 1   \\ 
 \hline
$-L(5,1)$    & 2  & 0 \\ 
 \hline
$-L(7,1)$    & 0  & 0 \\ 
\hline
\end{tabular}
\end{center}
The above $w(Y)$-values coincide with that of \cite{KMZ}. From (4), for example, in $-L(3,1)$ case, $w(X) + 3 \equiv \sigma(X)$ mod 4. $-L(p=odd,1)$ carry a unique spin structure $s_0$ and bound a smooth spin 4-manifold, hence, after applying the Rokhlin's theorem, we have $\mu (Y,s_0) = (w(X) + 3) $ mod 4. By the generalized signature theorem, we get $p_1 (TX) = 9	+ 3w(X) + 3\eta(0)$ mod 12. Hence, $X$ being a string manifold requires $w(X)=-3-\eta(0)$ mod 4. Otherwise, $9	+ 3w(X) + 3\eta(0) = 48\intg^{\ast}$~\cite{MK}. $d(Y)=1$ means that $Y$ possesses a 1-cocycle in $\intg_3$ coefficient group. We note that $d(Y_0)=2d(Y)$ and  $\text{def}_{3}(0)$ vanishes modulo 4. The latter implies that $w(\tilde{X})= 2w(X)$ since $F \cdot F = 0$ due to $\Theta =0$, which indicates that $Y$ has no effect on its bounding $X$ and similarly for the ($\tilde{X}, Y_0$) pair. We observe that $w(L(p,1))=1$ does not occur, which is also true when $p$ is even (see Appendix B). 
\newline

\noindent For $p=$ even, $H^{1}(Y;\intg/ 2\intg) = \intg/ 2\intg $ and hence we have
$$
\tau_{6}[Y, \Theta ] =  \sum_{b\in Spin^c(Y)/\intg_2} c^{Witt}_{\Theta b} \hat{Z}_{b} (q= e^{\frac{i2\pi}{6}}),\qquad \Theta = 0 , 1.
$$
Applying (9) and (12), we get
\begin{center}
\begin{tabular}{ |c|c| } 
\hline
$-L(p ,1 )$       & $\tau_{6}[\Theta=0] + \tau_{6}[\Theta=1]$ \\ 
 \hline
$-L(2,1)$ & $-1 + \sqrt{3}$  \\ 
 \hline
$-L(4,1)$ & $1- i\sqrt{3} $ \\ 
 \hline
$-L(6,1)$ & $i\sqrt{3} - \sqrt{3}  $ \\ 
 \hline
$-L(8,1)$ & $ -1 + i\sqrt{3}  $ \\ 
 \hline
$-L(10,1)$ & $ 1+ \sqrt{3} $ \\ 
 \hline
$-L(12,1)$ & $i\sqrt{3} - i\sqrt{3} $ \\ 
 \hline
$-L(14,1)$ & $-1 - \sqrt{3}$  \\ 
 \hline
$-L(16,1)$ & $1+i\sqrt{3}$  \\ 
 \hline
$-L(18,1)$ & $i\sqrt{3} + \sqrt{3}$  \\ 
 \hline
$-L(20,1)$ & $ -1 - i\sqrt{3}$  \\ 
 \hline
$-L(22,1)$ & $ 1 - \sqrt{3}$  \\ 
 \hline
$-L(24,1)$ & $i\sqrt{3} + i\sqrt{3} $  \\ 
\hline
\end{tabular}
\end{center}
This result is in agreement with that of \cite{KMZ}. The WRT invariant $\tau_{6}$ has a period from $p=2$ to $p=24$; in other words,
$$
\tau_{6}(p=26) = \tau_{6}(p=2),\quad \tau_{6}(p=28) = \tau_{6}(p=4)\quad \text{etc}.
$$
From the above table, $w(Y)$ and $\text{def}_{3}(\Theta)$ can be obtained, which are recorded in Appendix A. For $w(Y)=2$, we have $p_1 (TX)= 6 + 3w(X) + 3\eta(0)$ mod 12 by (4)  and the generalized signature theorem. For $w(X)=-2-\eta(0)$ mod 4, $X$ carries a string structure. The value of $def_3 (1)$ indicates how close or far are $\tilde{X}$ and $X$ being closed manifolds. When $def_3 (1)=0$, $Y$ and $Y_{1}$ have no effects on $X$ and $\tilde{X}$, respectively, since $w(\tilde{X})=2w(X)$ for genus zero $F$. This is the case for $-L(2,1)$ and $-L(22,1)$ whereas $-L(4,1)$ has an effect on its bounding $X$ and on $\tilde{X}$ via $Y_1$ since $w(\tilde{X}) \neq 2w(X)$ even for genus zero $F$. Therefore, $X$ bounded by $-L(4,1)$ and its cover $\tilde{X}$  are genuinely open manifolds. We cannot apply the Rokhlin's theorem to $p$ even cases since they may not bound a smooth 4-manifold.

\subsection{Other Seifert fibered manifolds}

Having analyzed special Seifert fibered manifolds in the previous section, we proceed with Seifert fibered manifolds with three or four singular fibers that are obtained from Dehn surgery on the trefoil or the figure eight knot. For the three singular fiber case, they can be obtained by an integer surgery on the right handed trefoil. 
$$
Y= S^{3}_{-2}(3^{r}_{1}) = M \lb -1 \bigg| \frac{1}{2}, \frac{1}{3}, \frac{1}{8} \rb, \qquad H_1(Y;\intg) = \intg/2\intg
$$
Its $\hat{Z}$ in terms of the Eichler integral of the weight 1/2 false theta functon are~\cite{CCFGH}
\begin{align*}
\hat{Z}_{0}[Y;q] & = q^{71/96} \lb \Psi_{24}^{(1)} - \Psi_{24}^{(17)} \rb  \\
\hat{Z}_{1}[Y;q] & =  -q^{71/96} \lb \Psi_{24}^{(7)} - \Psi_{24}^{(23)} \rb.
\end{align*}
Using (3), (9) and the formulas in Appendix B, we obtain
$$
\tau_6 [Y,0] = -1, \qquad \tau_6 [Y,1] = -3.
$$
From these values, we get via (5)
\begin{equation}
w(Y) = 2\quad \text{mod}\, 4, \qquad d(Y)=d(Y_0)=0, \qquad d(Y_1)=1,\qquad \text{def}_{3}(1) = 0\quad \text{mod}\, 4.
\end{equation}
The meaning of the first equation is explained in the $-L(p=\text{even},1)$ example. The third and the second equation indicate that there is a 1-cocycle in $Y_1$ generating $H^1 (Y_1;\intg/3\intg)$ whereas $Y$ and $Y_0$ have no such cocycles.  We observe that this manifold is the first instance in which $d(Y)=0$ but $d(Y_1) \neq 0$, which is absent in the case of Lens spaces. For $\text{def}_{3}(1) = 0$, $Y$ and $Y_1$ have no effect on how $\tilde{X}$ covers $X$. In other words, $\tilde{X}$ and $X$ behave as if they are closed.\\
\indent We can also compute $w(Y)$ using its 4-dimensional definition (4). In order to apply it, we need a plumbing graph $\Gamma$ description of $Y$, which is
\begin{center}
\begin{tikzpicture}
\centering
\tikzstyle{every node}=[draw,shape=circle]

\draw (0,0) node[circle,fill,inner sep=1pt,label=below:$-1$](-1){} -- (1,0) node[circle,fill,inner sep=1pt,label=right:$-8$](-8){};
\draw (0,0)  -- (0,1) node[circle,fill,inner sep=1pt,label=right:$-3$](-3){};
\draw (0,0) -- (-1,0) node[circle,fill,inner sep=1pt,label=left:$-2$](-2){};

\end{tikzpicture}
\end{center}
\noindent The adjacency matrix of $\Gamma$ is
$$
B(\Gamma) = \begin{pmatrix}
-1 & 1 & 1 & 1 \\
1 & -2 & 0 & 0 \\
1 & 0 & -3 & 0 \\
1 & 0 & 0 & -8 \\
\end{pmatrix}
$$
The compact oriented 4-manifold $X$ that is bounded by $Y$ is a negative definite graph 4-manifold that is characterized by $\Gamma$ as well. As a consequence, the intersection form of $X$ is given by $B(\Gamma)$. Upon application of (4), we arrive at the same value of $w(Y)$ as given in equation (13).
\newline 

\indent We next consider 
$$
Y= S^{3}_{-3}(3^{r}_{1}) = M \lb -1 \bigg| \frac{1}{2}, \frac{1}{3}, \frac{1}{9} \rb, \qquad H_1(Y;\intg) = \intg/3\intg
$$
Its $q$-series are~\cite{CCFGH}
\begin{align*}
\hat{Z}_{0}[Y;q] & = q^{71/72} \lb \Psi^{(1)}_{18} + \Psi^{(17)}_{18} \rb\\
\hat{Z}_{1}[Y;q] & = -q^{71/72} \lb \Psi^{(5)}_{18} + \Psi^{(13)}_{18} \rb.
\end{align*}
Using the values of the Eichler integrals of the $w=3/2$ false theta function at sixth primitive root of unity (see Appendix B) along with (9), we get
$$
\tau_6 [Y,0] = i \sqrt{3}.
$$
From (7), we arrive at
$$
w(Y)= 3\quad \text{mod}\, 4, \qquad d(Y) =1, \qquad d(Y_0)=2.
$$
The interpretation of the first equation is same as that of $-L(3,1)$ in Section 5.2. $d(Y) =1$ means that $Y$ has a 1-cocycle, which is lifted to two 1-cocycles in $Y_0$.\\
\indent We confirm the above $w(Y)$ using the same method as before together with
$$
B(\Gamma)= \begin{pmatrix}
-1 & 1 & 1 & 1 \\
1 & -2 & 0 & 0 \\
1 & 0 & -3 & 0 \\
1 & 0 & 0 & -9 \\
\end{pmatrix}
$$
\newline

Let us move to a manifold obtaind from Dehn surgery on the figure-eight knot $4_1$. In \cite{Thurston}, it was proved that the surgery on $4_1$ along an exceptional surgery slope produces a Seifert fibered rational homology sphere. We first choose 
$$
Y= S^{3}_{+3}(4_1) = M \lb -1 \bigg| \frac{1}{3}, \frac{1}{3}, \frac{1}{4} \rb
$$
We compute its $q$-series via (15) in Appendix C. To convert them in terms of the Eichler integrals of the false theta functions, we apply the method described in \cite{CCFGH}. Upon application, we find
\begin{align*}
\hat{Z}_0 & = q^{-1/48} \lb \Psi_{12,1}(q) - \Psi_{12,7}(q) \rb\\
\hat{Z}_1 & = -q^{-1/48} \Psi_{12,9}(q).
\end{align*}
Using (2) and (9) result in
$$
\tau_{6} [Y;0] = -i \sqrt{3}\quad \imp \quad d(Y)=1,\quad d(Y_0)=2,\quad w(Y)=1\quad \text{mod}\, 4.
$$
The last expression translates into $p_1 (X) = 3 + 3w(X)+ 3\eta(0)$ mod 12 by (4) and the generalized signature theorem, which means that for $w(X) \neq -1 -\eta(0)$ mod 4, $X$ is a not string manifold. In this case, $ 1 + w(X)+ \eta(0) \in 16\intg^{\ast}$ by the result of \cite{MK}. Interpretations of the other results are same as that of the previous example.
\newline

\noindent Another exceptional surgery is
$$
Y= S^{3}_{+2}(4_1) = M \lb -1 \bigg| \frac{1}{2}, \frac{1}{4}, \frac{1}{5} \rb
$$
Using the same procedure as the previous manifold yields
\begin{align*}
\hat{Z}_0 & =   q^{19/80} \lb \Psi_{12,1}(q) - \Psi_{12,9}(q) \rb\\
\hat{Z}_1 & =   - q^{19/80} \lb \Psi_{12,11}(q) - \Psi_{12,19}(q) \rb.
\end{align*}
From here, we obtain using (3), (9) and (7), 
$$
d(Y)=1,\quad d(Y_0)=2,\quad d(Y_1) = 0,\quad w(Y)=2\quad \text{mod}\, 4, \qquad \text{def}_{3}(1) = 2\quad \text{mod}\, 4.
$$
The first three values imply that $Y$ has one nontrivial 1-cocyle while $Y_0$ has two generators, which are a lift of the former cycle whereas $Y_1$ has none. The defect value essentially indicates that the deviation of $w(\tilde{X})$ from $2w(X)$ due to the presence of $Y$. Meaning of $w(Y)=2$ is explained in Section 5.2.
\newline

\noindent Another example of a hyperbolic knot is $5_2$. Performing $-3$ surgery on this knot produces
$$
Y= S^{3}_{-3}(5_2) = M \lb -2 \bigg| \frac{2}{3}, \frac{2}{3}, \frac{3}{5} \rb
$$
Its plumbing graph is
\begin{center}
\begin{tikzpicture}
\centering
\tikzstyle{every node}=[draw,shape=circle]

\draw (0,0) node[circle,fill,inner sep=1pt,label=below:$-2$](-2){} -- (1,0) node[circle,fill,inner sep=1pt,label=below:$-2$](-2){}-- (2,0) node[circle,fill,inner sep=1pt,label=below:$-2$](-2){};
\draw (0,0)  -- (0,1) node[circle,fill,inner sep=1pt,label=right:$-2$](-2){}-- (0,2) node[circle,fill,inner sep=1pt,label=right:$-2$](-2){};
\draw (0,0) -- (-1,0) node[circle,fill,inner sep=1pt,label=below:$-2$](-2){}-- (-2,0) node[circle,fill,inner sep=1pt,label=below:$-3$](-3){};

\end{tikzpicture}
\end{center}
Through the same method, we get
$$
\hat{Z}_0 =   q^{-16/15} \lb \Psi_{15,2}(q) - \Psi_{15,8}(q) \rb, \qquad \hat{Z}_1  =   - q^{-16/15}  \Psi_{15,12}(q)\\
$$
and 
$$
d(Y)=1,\quad d(Y_0)=2,\quad w(Y)=3\quad \text{mod}\, 4.
$$
This result is same as that of the second example in this section.
\newline

We finish this section with a Seifert manifold with four singular fibers. This manifold has a mixed modular property compared with a three singular fiber case. Specifically, in addition to the presence of the Eichler integral of the weight $3/2$ false theta function, the weight $1/2$ false theta function appears as well~\cite{H1}. We consider 
$$
Y = M \lb -2 \bigg| \frac{1}{2}, \frac{2}{3}, \frac{2}{5}, \frac{2}{5} \rb \qquad H_{1}(Y,\intg)= \intg / 5 \intg 
$$
Its $\hat{Z}$ are given in \cite{CCFGH}\footnote{This a typo in Section 8 of \cite{CCFGH}. I would like to thank Sarah Harrison for informing the correct expression.} 
\begin{align*}
\hat{Z}_{0}(q) & = q^{-109/120} \lb \Psi_{30,23} - B_{30,7} + B_{30,13} - B_{30,17} + B_{30,23} \rb \\
\hat{Z}_{1}(q) & = 0\\
\hat{Z}_{2}(q) & = 2 q^{-109/120} \lb  B_{30,5} - B_{30,25} \rb\\
B_{m,r} (q)   & := \frac{1}{2m} \big[ \Phi_{m,r}(q) - r \Psi_{m,r}(q) \big],
\end{align*}
where $\Phi_{m,r}(q)$ and $\Psi_{m,r}(q)$ are the Eichler integral of the weight $1/2$ and $3/2$  false theta functions, respectively.
After using (7) and the formulas in Appendix B, we find
$$
w(Y) = 2\quad \text{mod}\, 4 \qquad d(Y) = 0.
$$
They coincide with the first two results in (13), hence their meaning follows from there.

\subsection{Hyperbolic manifold}

\begin{figure}[t]
\subfloat{\includegraphics[scale = 0.28]{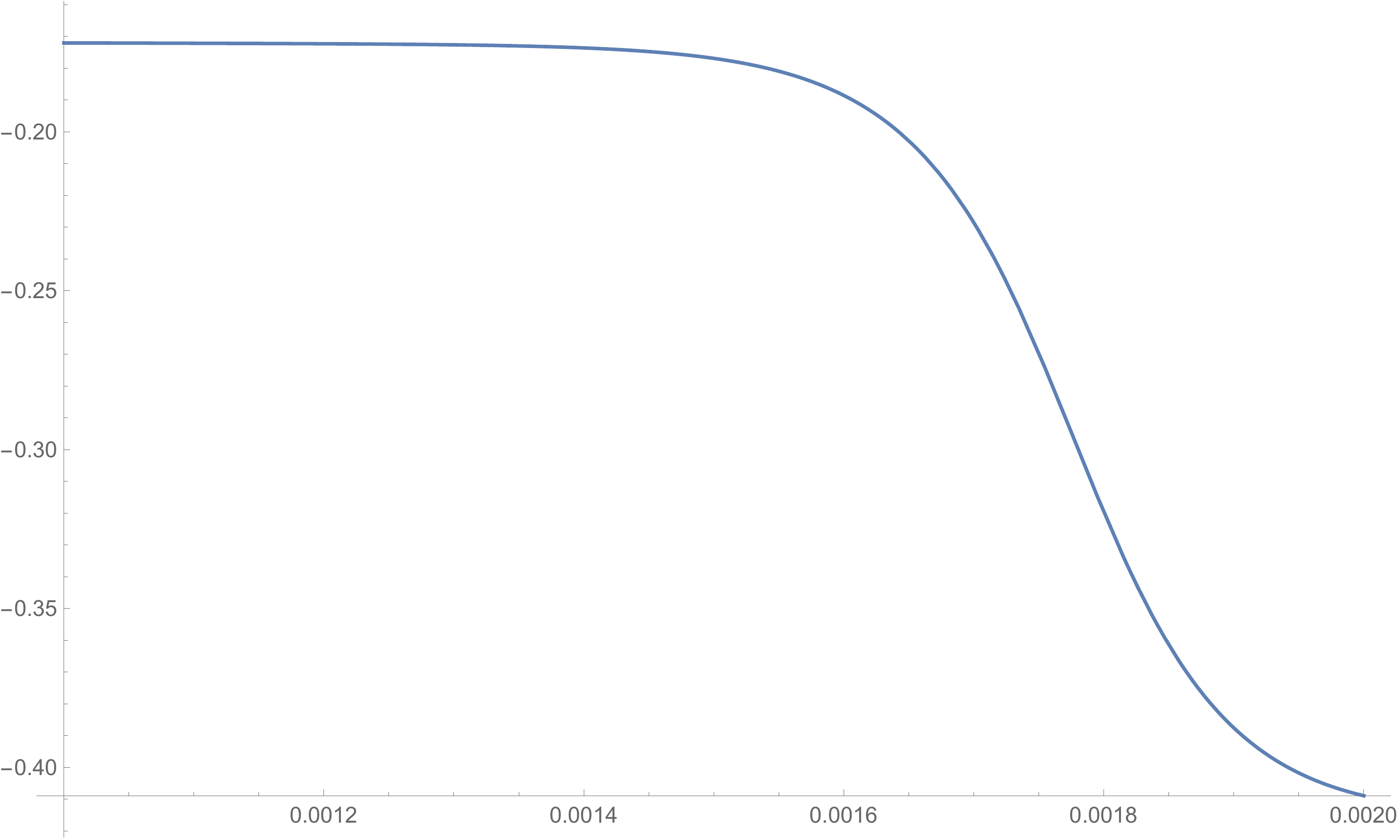}}
\hfill
\subfloat{\includegraphics[scale = 0.2]{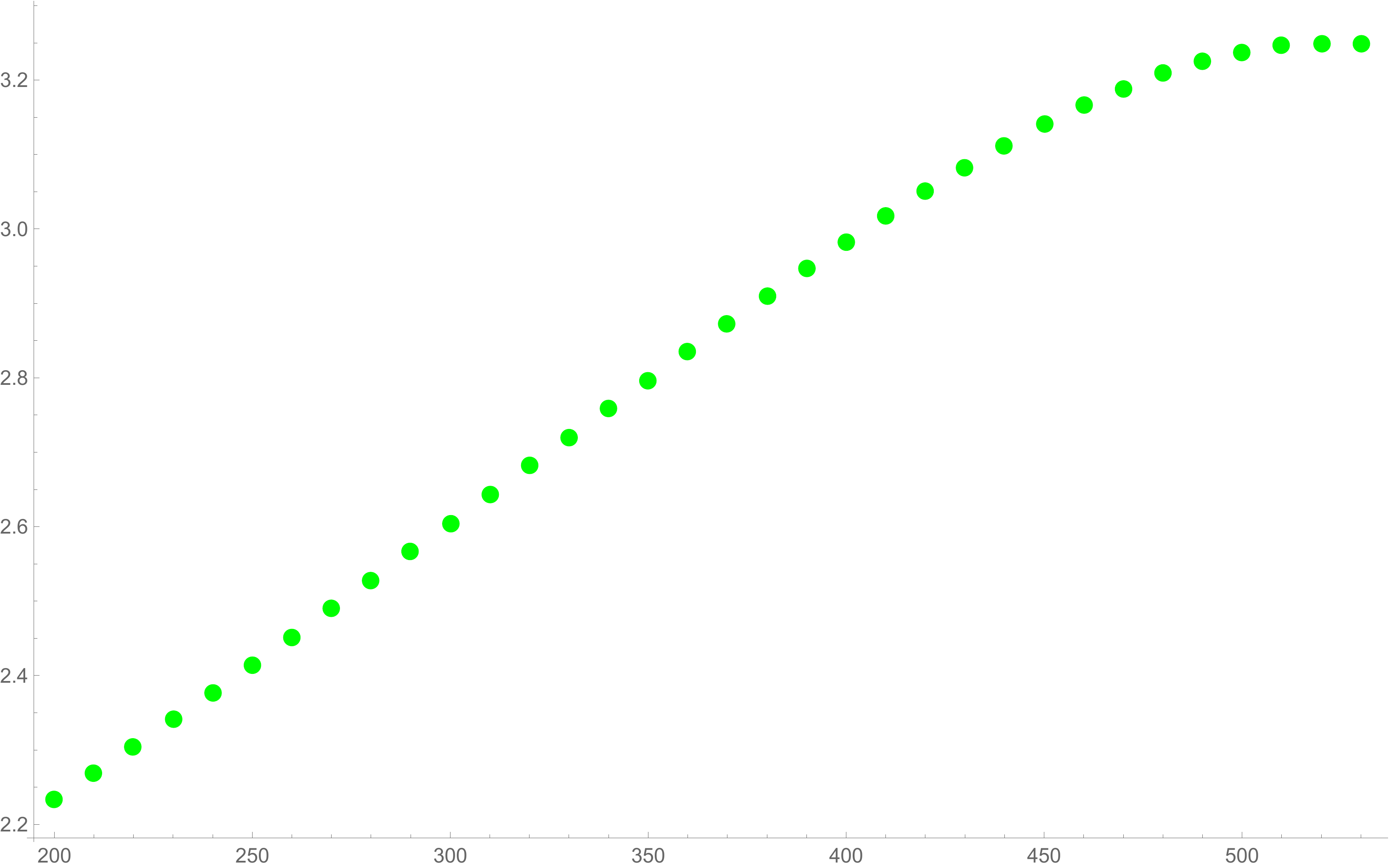}}
\caption{The limiting behaviors of the argument of $-q^{1/2} \hat{Z}[S^{3}_{-1/2}(4_1)] $ (Left) and the absolute value (Right) of (14) as q goes to the sixth root of unity. The left plot corresponds to a truncation at $1800$ and its horizontal axis is the y-range, $0.001\leq y \leq 0.002$. The horizontal axis of the right plot is the truncation powers of $q$. The normalized argument and the absolute value approach $-0.17$ and $3.248$, respectively.}
\end{figure}

In three dimensions, hyperbolic knots and manifolds are abundant. The latter can be obtained using the former and the Dehn surgery performed along a generic slope~\cite{Th2}. For a hyperbolic manifold of interest, we consider the $-1/2$ surgery on the figure eight knot $K=4_1$. Since this surgery slope is outside the set of the exceptional slopes for $K$, the resulting manifold possess a hyperbolic structure. Its $\hat{Z}$ can be easily found using the recent results in knot complement series in \cite{GM}.
\begin{align}
\hat{Z}[S^{3}_{-\frac{1}{2}}(4_1);q] & = -q^{-\frac{1}{2}} \lb 1-q+2 q^3-2 q^6+q^9+3 q^{10}+q^{11}-q^{14}-3 q^{15}-q^{16}+2 q^{19}+2 q^{20} + \cdots  \notag\right.\\
&   +2335418615 q^{1600} + \cdots  \left.\rb.
\end{align} 
The nontriviality is finding a limiting value of $\hat{Z}$ of a hyperbolic manifold as $q$ goes to a root of unity. The difficulty comes from rapidly growing coefficients of (14). Although a closed form formula for $\hat{Z}$ of $-1/r$-surgery on any knot is available in \cite{EGGKPSS}, finding the limit poses its own challenges. We leave it for future work and use a numerical method, which was employed in \cite{GPP}. In this method, we approximate (14) by a finite number of terms and we set $q=e^{i 2\pi (\frac{1}{6} + iy)},\, y\in \real_{+}$, where $q$ is in an unit disk in the complex plane (see Figure 1). We then analyze the radial limit $y \rarw 0$ behavior of (14). We iterate this procedure by varying the number of truncated terms to find the limiting value of (14). 
\newline

\noindent The argument part of (14) is
\begin{equation}
\frac{1}{2\pi} \text{Arg}  \lb \hat{Z}[S^{3}_{-\frac{1}{2}}(4_1)] \rb = \frac{5}{12} + \frac{1}{2\pi} \text{Arg} \lb -q^{1/2} \hat{Z}[S^{3}_{-\frac{1}{2}}(4_1)] \rb,
\end{equation}
where $5/12$ is from $-q^{-1/2}$. For the second argument, several truncations of (14) were analyzed for $0.001 \leq y \leq 0.002$. Plots having similar behavior to Figure 3 (left) also occur at other truncation powers, for instance, $1860, 1980$, and $2040$. They all approach $-0.17$ as $y$ goes to the lower bound. For the absolute value of (14), truncations of the powers between $200$ and $530$ are considered, the absolute value approaches $3.248$ (Figure 3). Using the limiting values, we get
\begin{equation}
\tau_6 [S^{3}_{-\frac{1}{2}}(4_1)] = \frac{-i}{2} \lim_{q \rarw e^{\frac{i 2\pi}{6}}}  \hat{Z}[S^{3}_{-1/2}(4_1);q] \approx 1.62 .
\end{equation}
According to \cite{KMZ}, $\tau_6 [\intg HS]=1$. Although the absolute value approximation may look crude, the coefficients of (14) increases exponentially fast, for example, it is order of $10^{10}$ at $q^{2040}$, however, (16) is $O(1)$. It turns out (16) being a real number, which is solely determined by the phase value (15) and the surgered manifold being $\intg HS$ are sufficient to arrive at 
$$
w(Y) = 0\quad \text{mod}\, 4 \qquad d(Y) = 0, \qquad d(Y_0) = 0.
$$
Their interpretations can be found in Section 5.1.



\section{Open Questions}

We list open questions.

\begin{itemize}

	\item Evaluation of a limit of $\hat{Z}$ at a root of unity is nontrivial in general, which is required for finding Witt invariants, Rokhlin invariant~\cite{GPP} and WRT invariant~\cite{GPPV}. For Seifert manifolds with three singular fibers, a method for finding their $\hat{Z}$'s in terms of the false theta functions is available~\cite{CCFGH}. Hence, finding a limit of $\hat{Z}$ at a root of unity can be done analytically. However, for Seifert manifolds with more than three singular fibers, a method for finding their $\hat{Z}$'s in terms of the false theta functions has not yet been found.
	
	\item  Although a closed form formula for $\hat{Z}$ of hyperbolic 3-manifolds exist as mentioned in Section 5.4, analytic evaluation of the limit is nontrivial. Finding an analytic method would enlarge the range of 3-manifolds whose Witt invariants, Rokhlin invariant and correction term (d-invariant) can be found. Furthermore, computations of the invariants for hyperbolic 3-manifolds that are rational homology sphere have not been explored.   
	
\item Other invariants of 3-manifolds at different roots of unity are mentioned in \cite{KM}. It would be interesting to find connections between them and $\hat{Z}$.

\end{itemize}

\noindent \textbf{Acknowledgments.} I would like to thank Sungbong Chun, Sarah Harrison, Kazuhiro Hikami, Robion Kirby, Paul Melvin and Pavel Putrov for helpful explanations. I am grateful to Sergei Gukov for numerous explanations and the suggestion on this manuscript. I would also like to thank the referee for the suggestions that led to an improvement of my manuscript. 

\appendix
\section*{Appendix}
\addcontentsline{toc}{section}{Appendix}

\section{Witt invariants for the Lens spaces}

We summarize the Witt invariants for  $L(p,1)$, where $p$ is even. 
\begin{center}
\begin{tabular}{ |c|c|c|c|c|c| } 
\hline
$-L(p ,1 )$  & $w(Y) \in \intg /4\intg $ & $d(Y) \in \intg$  & $d(Y_1) \in \intg$ & $\text{def}_{3}(1)  \in \intg /4\intg $ & $ 2(1^3) \in \intg /4\intg $   \\ 
\hline
$-L(2,1)$  &  $2$  & $0$ & $0$ & $0$ & $2$ \\ 
\hline
$-L(4,1)$  &  0  & 0 & 0 & 3 & 0\\ 
\hline
$-L(6,1)$  &  3  & 1 & 1 & 3 & 2 \\ 
\hline
$-L(8,1)$  &  2  & 0 & 0 & 3 & 0 \\
\hline 
$-L(10,1)$ & 0 & 0 & 0 & 2 & 2 \\ 
\hline
$-L(12,1)$ & 3 & 1 & 1 & 2 & 0 \\ 
\hline
$-L(14,1)$ & 2 & 0 & 0 & 2 & 2  \\ 
\hline
$-L(16,1)$ & 0 & 0 & 0 & 1 & 0  \\ 
\hline
$-L(18,1)$ & 3 & 1 & 1 & 3 & 0  \\
\hline 
$-L(20,1)$ & 2 & 0 & 0 & 1 & 0  \\ 
\hline
$-L(22,1)$ & 0 & 0 & 0 & 0 & 2  \\ 
\hline
$-L(24,1)$ & 3 & 1 & 1 & 2 & 2  \\ 
\hline
\end{tabular}
\end{center}
As written in Section 5.2, $d(Y_0) = 2d(Y)$ and $\text{def}_{3}(0) = 0$ modulo 4. 

\section{$\hat{Z}$-series for Brieskorn spheres and modular forms}

We record $\hat{Z}$ for the manifolds in (12) and the formulas for the weight $1/2$ and $3/2$ modular forms at k-th root of unity.

$$
\hat{Z} [S^{3}_{-\frac{1}{r}}(T(2,5)) ;q] = q^{\frac{71}{40}-\frac{r}{4}-\frac{5}{2+20 r}} \lb \Psi^{(30r -7)}_{100r+10} (q) - \Psi^{(30r +13 )}_{100r+10}(q) - \Psi^{(70r -3)}_{100r+10}(q) + \Psi^{(70r +17)}_{100r+10}(q) \rb
$$

$$
\hat{Z} [S^{3}_{-\frac{1}{r}}(T(2,7)) ;q] = q^{\frac{143}{56}-\frac{r}{4}-\frac{7}{2+28 r}} \lb \Psi^{(70r -9)}_{196r+14} (q) - \Psi^{(70r +19)}_{196r+14}(q) - \Psi^{(126r -5)}_{196r+14}(q) + \Psi^{(126r +23)}_{196r+14}(q) \rb
$$

$$
\hat{Z} [S^{3}_{-\frac{1}{r}}(T(3,5)) ;q] = q^{\frac{191}{60}-\frac{r}{4}-\frac{15}{4+60 r}} \lb \Psi^{(105r -8)}_{225r+15} (q) - \Psi^{(105r +22)}_{225r+15}(q) - \Psi^{(195 r-2 )}_{225r+15}(q) + \Psi^{(195 r + 28)}_{225r+15}(q) \rb
$$

$$
\hat{Z} [S^{3}_{-\frac{1}{r}}(T(3,7)) ;q] = q^{\frac{383}{84}-\frac{r}{4}-\frac{21}{4+84 r}} \lb \Psi^{(231 r -10)}_{441r+21} (q) - \Psi^{(231 r +32)}_{441r+21}(q) - \Psi^{(357 r-4 )}_{441r+21}(q) + \Psi^{(357 r + 38)}_{441r+21}(q) \rb
$$

$$
\hat{Z} [S^{3}_{-\frac{1}{r}}(T(2,-5)) ;q] = q^{-\frac{71}{40}+\frac{5}{2-20 r}-\frac{r}{4}} \lb \Psi^{(30r -13)}_{100r-10} (q) - \Psi^{(30r +7 )}_{100r-10}(q) - \Psi^{(70r -17)}_{100r-10}(q) + \Psi^{(70r +3)}_{100r-10}(q) \rb
$$

$$
\hat{Z} [S^{3}_{-\frac{1}{r}}(T(2,-7)) ;q] = q^{-\frac{143}{56}+\frac{7}{2-28 r}-\frac{r}{4}} \lb \Psi^{(70r -19)}_{196r-14} (q) - \Psi^{(70r + 9)}_{196r-14}(q) - \Psi^{(126r -23)}_{196r-14}(q) + \Psi^{(126r +5)}_{196r-14}(q) \rb
$$

$$
\hat{Z} [S^{3}_{-\frac{1}{r}}(T(3,-5)) ;q] = q^{-\frac{191}{60}+\frac{15}{4-60 r}-\frac{r}{4}} \lb \Psi^{(105r -22)}_{225r-15} (q) - \Psi^{(105r +8)}_{225r-15}(q) - \Psi^{(195 r-28 )}_{225r-15}(q) + \Psi^{(195 r + 2)}_{225r-15}(q) \rb
$$

$$
\hat{Z} [S^{3}_{-\frac{1}{r}}(T(3,-7)) ;q] = q^{-\frac{383}{84}+\frac{21}{4-84 r}-\frac{r}{4}} \lb \Psi^{(231 r -32)}_{441r-21} (q) - \Psi^{(231 r +10)}_{441r-21}(q) - \Psi^{(357 r- 38 )}_{441r-21}(q) + \Psi^{(357 r + 4)}_{441r-21}(q) \rb
$$
\newline

The Eichler integrals of the $w=1/2$ and $w=3/2$ false theta functions $\Phi_{m,r}(q)$ and  $\Psi_{m,r}(q)$ at k-th primitive root of unity are given by~\cite{H1, H2}.
\begin{align*}
\Phi_{m,r}(e^{i 2 \pi /k}) & = -(m k) \sum _{n=1}^{2 m k} \left(\left(\frac{n}{2 m k}\right)^2-\frac{n}{2 m k}+\frac{1}{6}\right) \psi^{\prime, (r)}_{2m}(n) e^{\frac{i \pi  n^2}{2 m k}}\\
\psi^{\prime, (r)}_{2m}(n) & : = \begin{cases}
 1, n \equiv \pm r\quad \text{mod}\quad 2m\\
0, \text{otherwise}
\end{cases}
\end{align*}

\begin{align*}
\Psi^{(r)}_{m} \equiv \Psi_{m,r}(e^{i 2 \pi /k}) & = \sum _{n=1}^{2 m k} \left(  \frac{1}{2} - \frac{n}{2mk} \right) \psi^{(r)}_{2m}(n) e^{\frac{i \pi  n^2}{2 m k}} \\
\psi^{(r)}_{2m}(n) & : = \begin{cases}
\pm 1, n \equiv \pm r\quad \text{mod}\quad 2m\\
0, \text{otherwise}
\end{cases}
\end{align*}

\section{$\hat{Z}$-series for plumbed 3-manifolds}

We state a formula for $\hat{Z}$ of (weakly) negative definite plumbed 3-manifolds $Y(\Gamma)$ having $b_{1} (Y(\Gamma))=0$~\cite{GPPV, GM}:
\begin{equation}
\hat{Z}_{b}[Y(\Gamma) ;q]= (-1)^{\pi} q^{\frac{3\sigma - Tr B}{4}}\, \prod_{v \in Vert} PV \oint_{|z_v|=1}  \frac{d z_v}{i2\pi z_v} \lb z_v - \frac{1}{z_v} \rb^{2-deg(v)} \Theta^{Y}_{b}(\vec{z};q),
\end{equation}
where
$$
\Theta^{Y}_{b}(\vec{z};q) = \sum_{\vec{w} \in 2B\intg^{L} + \vec{b}} q^{-\frac{(\vec{w},B^{-1}\vec{w})}{4}} \prod_{v \in Vert} z_{v}^{w_v},\qquad b \in Spin^c(Y)\cong H_1(Y)
$$
$$
B= \text{adjacency matrix of}\, \Gamma,\qquad    \pi = \text{$\sharp$ (positive eigenvalues of B)},\qquad \sigma = \text{signature}(B),
$$
$$
PV= \lim_{\epsilon \rarw 0 } \frac{1}{2} \lb \oint_{|z_v|=1 + \epsilon} + \oint_{|z_v|=1 - \epsilon} \rb
$$

\end{document}